\documentclass[11pt, twoside]{article}
\usepackage{latexsym}
\usepackage{amsmath}\def\rad{\mbox{rad}}
\usepackage{amssymb}
\usepackage[all]{xy}
\usepackage{amsfonts}
\usepackage{verbatim}
\usepackage{amsthm}
\usepackage{mathrsfs}
\usepackage{epsfig}
\usepackage{xy}
\usepackage{array}
\usepackage{stmaryrd}
\usepackage{graphicx,color}
\usepackage{xcolor}
\usepackage[colorlinks=true,linkcolor=blue,citecolor=blue]{hyperref}
\usepackage{tikz}
\usetikzlibrary{arrows,calc}
\usepackage{etex}
\usepackage{mathdots}
\usepackage{float}
\usepackage{graphics}
\usepackage{pdflscape}
\usepackage{bm}
\usepackage{anysize,hyperref}
\input xypic
\xyoption{all}
\newcommand{\mr}{\hbox{\boldmath$\cdot$}}
\usepackage[perpage,symbol]{footmisc}
\usepackage{setspace}
\def\A{\mathscr{A}}

\def\P{\mathcal{P}}
\def\I{\mathcal{I}}
\def\B{\mathscr{B}}
\def\C{\mathscr{C}}
\def\E{\mathbb{E}}
\def\F{\mathbb{F}}
\def\s{\mathfrak{s}}\def\t{\mathfrak{t}}
\def\id{\mathrm{id}}
\def\op{^\mathrm{op}}
\def\Ab{\mathit{Ab}}
\def\del{\delta}
\def\dr{\ar@{->}[r]}

\def\X{\mathscr{X}}
\def\Cone{\mbox{Cone}}
\def\add{\mbox{add}}
\def\CoCone{\mbox{CoCone}}
\def\Aut{\mbox{Aut}}

\newcommand{\CC}{{\bf{C}}^{n+2}_{\C}}
\def\add{\mathsf{add}\hspace{.01in}}
\newcommand{\ov}{\overset}
\newcommand{\lra}{\longrightarrow}
\newcommand{\co}{\colon}
\newcommand{\uas}{^{\ast}}            %%% ^*
\newcommand{\sas}{_{\ast}}
\newcommand{\Xd}{\langle X_{\bullet},\del\rangle}  %%% <X,¦Ä>
\newcommand{\Yr}{\langle Y_{\bullet},\rho\rangle}  %%% <Y,¦Ñ>
\newcommand{\ush}{^\sharp}           %%% ^sharp
\newcommand{\ssh}{_\sharp}
\SelectTips{cm}{10}

\begin{document}
\title{\Large{\bf $\bm{n}$-exact categories arising from $\bm{n}$-exangulated categories\footnotetext{\hspace{-1em} Panyue Zhou was supported by the National Natural Science Foundation of China (Grant No. 11901190) and the Scientific Research Fund of Hunan Provincial Education Department (Grant No. 19B239).}}}
\medskip
\author{Jian He and Panyue Zhou}

\date{}

\maketitle
\def\blue{\color{blue}}
\def\red{\color{red}}

\newtheorem{theorem}{Theorem}[section]
\newtheorem{lemma}[theorem]{Lemma}
\newtheorem{corollary}[theorem]{Corollary}
\newtheorem{proposition}[theorem]{Proposition}
\newtheorem{conjecture}{Conjecture}
\theoremstyle{definition}
\newtheorem{definition}[theorem]{Definition}
\newtheorem{question}[theorem]{Question}
\newtheorem{remark}[theorem]{Remark}
\newtheorem{remark*}[]{Remark}
\newtheorem{example}[theorem]{Example}
\newtheorem{example*}[]{Example}
\newtheorem{condition}[theorem]{Condition}
\newtheorem{condition*}[]{Condition}
\newtheorem{construction}[theorem]{Construction}
\newtheorem{construction*}[]{Construction}

\newtheorem{assumption}[theorem]{Assumption}
\newtheorem{assumption*}[]{Assumption}

\baselineskip=17pt
\parindent=0.5cm
\vspace{-6mm}

\begin{abstract}
\baselineskip=16pt
Let $\C$ be a Krull-Schmidt $n$-exangulated category and $\A$ be an $n$-extension closed subcategory of $\C$. Then $\A$ inherits the $n$-exangulated structure from the given $n$-exangulated category in a natural way. This construction gives $n$-exangulated categories which are not $n$-exact categories in the sense of Jasso nor $(n+2)$-angulated categories  in the sense of Geiss--Keller--Oppermann in general. Furthermore, we also give a sufficient
condition on when an $n$-exangulated category $\A$ is an $n$-exact category.
 These results generalize work by Klapproth and Zhou.\\[0.2cm]
\end{abstract}

\pagestyle{myheadings}
\markboth{\rightline {\scriptsize   J. He and P. Zhou}}
         {\leftline{\scriptsize $n$-exact categories arising from $n$-exangulated categories}}

\section{Introduction}
The notion of extriangulated categories was introduced by Nakaoka and Palu in \cite{NP} as a simultaneous generalization of
exact categories and triangulated categories. Exact categories (abelian categories are also exact categories) and triangulated categories are extriangulated categories, while there are some other examples of extriangulated categories which are neither exact nor triangulated, see \cite{NP,ZZ,HZZ,NP1}. Nakaoka and Palu \cite[Remark 2.18]{NP} proved extension closed subcategories of extriangulated categories are extriangulated categories. For example, let $A$ be an artin algebra and $K^{[-1,\hspace{0.5mm}0]}({\rm proj}A)$ the category of complexes of finitely  generated projective $A$-modules concentrated in degrees $-1$ and $0$, with morphisms considered up to homotopy. Then $K^{[-1,\hspace{0.5mm}0]}({\rm proj}A)$  is an extension closed subcategory of the bounded homotopy category $K^b({\rm proj}A)$ which is not exact and  triangulated, see \cite[Example 6.2]{INP}.
This construction gives extriangulated categories which are not exact and triangulated.

In \cite{GKO}, Geiss, Keller and Oppermann introduced
a new type of categories, called $(n+2)$-angulated categories, which generalize
triangulated categories: the classical triangulated categories are the special case $n=1$. They appear for example as $n$-cluster tilting subcategories of triangulated categories which are closed under the $n$th power of the shift functor.
Later, Jasso \cite{Ja} introduced $n$-exact categories which are categories inhabited
by certain exact sequences with $n+2$ terms, called $n$-exact sequences. The case $n=1$
corresponds to the usual concepts of exact categories.
An important source of examples of $n$-exact categories are $n$-cluster
tilting subcategories, see \cite[Theorem 4.14]{Ja}.
Recently, Herschend, Liu and Nakaoka \cite{HLN} introduced the notion of $n$-exangulated categories. It should be noted that the case $n =1$ corresponds to extriangulated categories. As typical examples we have that $n$-exact ($n$-abelian categories are also $n$-exact categories) and $(n+2)$-angulated categories are $n$-exangulated, see \cite[Proposition 4.5 and Proposition 4.34]{HLN}.  There are also some $n$-exangulated categories which are neither $n$-exact nor
 $(n+ 2)$-angulated, see \cite{HLN,LZ,HZZ1}.

Let $(\C,\E,\s)$ be an $n$-exangulated category and $\textbf{C}_{\C}$ denote the category of complexes in $\C$. We define $\textbf{C}^{n+2}_{\C}$ to be the full
subcategory of $\textbf{C}_{\C}$ consisting of complexes whose components are zero in all degrees outside of $\{0,1,\cdots,n+1\}$.
Let $\A$ be an $n$-extension closed subcategory of $\C$. We define $\E_{\A}$ to be the restriction of $\E$ onto $\A^{\rm op}\times\A$.
For any $\delta\in\E_{\A}(C,A)$, take a distinguished $n$-exangle $\langle X_{\bullet},\del \rangle$ with $X_i$ in $\A$ for $i=1,\cdots,n$. Define $\t(\delta)=[X_{\bullet}]$, where the homotopy  equivalence class is taken in ${\bf{C}}^{n+2}_{({\A;\hspace{0.5mm} A,\hspace{0.5mm} C})}$.

Herschend, Liu and Nakaoka \cite[Proposition 2.35]{HLN} proved that if $\t$-inflations are closed under composition and $\t$-deflations are closed under composition, then $(\A,\E_{\A},\t)$
is an $n$-exangulated category. Zhou \cite[Theorem 3.4]{Z} recently also proved that
this result holds in an $(n+2)$-angulated category, but he showed that
 this hypothesis ($\t$-inflations are closed under composition and $\t$-deflations are closed under composition) of Herschend, Liu and Nakaoka \cite[Proposition 2.35]{HLN} is automatically satisfied
 in an $(n+2)$-angulated category.
Based on this idea, we prove the first main result in the article, which is a higher couterpart of Nakaoka-Palu's result. Moreover, the construction gives $n$-exangulated categories which are not
$n$-exact nor $(n+2)$-angulated in general.
\begin{theorem}{\rm (see Theorem \ref{main1} for details)}\label{main111}
Let $(\C,\E,\s)$ be a Krull-Schmidt $n$-exangulated categor and
$\A$ be an $n$-extension closed subcategory of $\C$. Then $(\A,\E_{\A},\t)$
is an $n$-exangulated category.
\end{theorem}

Let $(\C,\Sigma,\Theta)$ be a  Krull-Schmidt $(n+2)$-angulated category and $\A$ be an $n$-extension closed subcategory of $\C$. Klapproth \cite{K} proved that $\A$ has the structure of an $n$-exact category under a natural assumption. Herschend, Liu and Nakaoka \cite[Proposition 4.37]{HLN} gave a description of when an $n$-exangulated category can become $n$-exact category. Based on this fact, combined with Theorem \ref{main111}, we prove the second main result in the article. Moreover, the construction gives a kind of example of $n$-exact categories.

\begin{theorem}{\rm (see Theorem \ref{main2} for details)}\label{main11}
Let $(\C,\E,\s)$ be a Krull-Schmidt $n$-exangulated category with enough projectives and enough injectives,
$\A$ be an $n$-extension closed subcategory of $\C$.
If $\C(\Sigma\A,\A)=\C(\A,\Omega\A)=0$, then $(\A,\mathscr{E}_{\A})$
is an $n$-exact category.
\end{theorem}

 This extends a recent main result of  Klapproth \cite{K} for $(n+2)$-angulated categories. Moreover, our proof method is to avoid proving that $\mathscr{E}_{\A}$ is closed under weak isomorphisms, which is very complicated and difficult.

This article is organized as follows. In Section 2, we review some elementary definitions and facts on $n$-exangulated categories. In Section 3, we prove our first and second main results.

\section{Preliminaries}
In this section, let $\C$ be an additive category and $n$ be a positive integer. Suppose that $\C$ is equipped with an additive bifunctor $\E\colon\C\op\times\C\to{\rm Ab}$, where ${\rm Ab}$ is the category of abelian groups. Next we briefly recall some definitions and basic properties of $n$-exangulated categories from \cite{HLN}. We omit some
details here, but the reader can find them in \cite{HLN}.

{ For any pair of objects $A,C\in\C$, an element $\del\in\E(C,A)$ is called an {\it $\E$-extension} or simply an {\it extension}. We also write such $\del$ as ${}_A\del_C$ when we indicate $A$ and $C$. The zero element ${}_A0_C=0\in\E(C,A)$ is called the {\it split $\E$-extension}. For any pair of $\E$-extensions ${}_A\del_C$ and ${}_{A'}\del{'}_{C'}$, let $\delta\oplus \delta'\in\mathbb{E}(C\oplus C', A\oplus A')$ be the
element corresponding to $(\delta,0,0,{\delta}{'})$ through the natural isomorphism $\mathbb{E}(C\oplus C', A\oplus A')\simeq\mathbb{E}(C, A)\oplus\mathbb{E}(C, A')
\oplus\mathbb{E}(C', A)\oplus\mathbb{E}(C', A')$.

For any $a\in\C(A,A')$ and $c\in\C(C',C)$,  $\E(C,a)(\del)\in\E(C,A')\ \ \text{and}\ \ \E(c,A)(\del)\in\E(C',A)$ are simply denoted by $a_{\ast}\del$ and $c^{\ast}\del$, respectively.

Let ${}_A\del_C$ and ${}_{A'}\del{'}_{C'}$ be any pair of $\E$-extensions. A {\it morphism} $(a,c)\colon\del\to{\delta}{'}$ of extensions is a pair of morphisms $a\in\C(A,A')$ and $c\in\C(C,C')$ in $\C$, satisfying the equality
$a_{\ast}\del=c^{\ast}{\delta}{'}$.}
Then the functoriality of $\E$ implies $\E(c,a)=a_{\ast}(c^{\ast}\del)=c^{\ast}(a_{\ast}\del)$.

\begin{definition}\cite[Definition 2.7]{HLN}
Let $\bf{C}_{\C}$ be the category of complexes in $\C$. As its full subcategory, define $\CC$ to be the category of complexes in $\C$ whose components are zero in the degrees outside of $\{0,1,\ldots,n+1\}$. Namely, an object in $\CC$ is a complex $X_{\bullet}=\{X_i,d_i^X\}$ of the form
\[ X_0\xrightarrow{d_0^X}X_1\xrightarrow{d_1^X}\cdots\xrightarrow{d_{n-1}^X}X_n\xrightarrow{d_n^X}X_{n+1}. \]
We write a morphism $f_{\bullet}\co X_{\bullet}\to Y_{\bullet}$ simply $f_{\bullet}=(f_0,f_1,\ldots,f_{n+1})$, only indicating the terms of degrees $0,\ldots,n+1$.
\end{definition}

\begin{definition}\cite[Definition 2.11]{HLN}
By Yoneda lemma, any extension $\del\in\E(C,A)$ induces natural transformations
\[ \del\ssh\colon\C(-,C)\Rightarrow\E(-,A)\ \ \text{and}\ \ \del\ush\colon\C(A,-)\Rightarrow\E(C,-). \]
For any $X\in\C$, these $(\del\ssh)_X$ and $\del\ush_X$ are given as follows.
\begin{enumerate}
\item[\rm(1)] $(\del\ssh)_X\colon\C(X,C)\to\E(X,A)\ :\ f\mapsto f\uas\del$.
\item[\rm (2)] $\del\ush_X\colon\C(A,X)\to\E(C,X)\ :\ g\mapsto g\sas\delta$.
\end{enumerate}
We simply denote $(\del\ssh)_X(f)$ and $\del\ush_X(g)$ by $\del\ssh(f)$ and $\del\ush(g)$, respectively.
\end{definition}

\begin{definition}\cite[Definition 2.9]{HLN}
 Let $\C,\E,n$ be as before. Define a category $\AE:=\AE^{n+2}_{(\C,\E)}$ as follows.
\begin{enumerate}
\item[\rm(1)]  A pair $\Xd$ is an object of the category $\AE$ with $X_{\bullet}\in\CC$
and $\del\in\E(X_{n+1},X_0)$, called an $\E$-attached
complex of length $n+2$, if it satisfies
$$(d_0^X)_{\ast}\del=0~~\textrm{and}~~(d^X_n)^{\ast}\del=0.$$
We also denote it by
$$X_0\xrightarrow{d_0^X}X_1\xrightarrow{d_1^X}\cdots\xrightarrow{d_{n-2}^X}X_{n-1}
\xrightarrow{d_{n-1}^X}X_n\xrightarrow{d_n^X}X_{n+1}\overset{\delta}{\dashrightarrow}.$$
\item[\rm (2)]  For such pairs $\Xd$ and $\langle Y_{\bullet},\rho\rangle$,  $f_{\bullet}\colon\Xd\to\langle Y_{\bullet},\rho\rangle$ is
defined to be a morphism in $\AE$ if it satisfies $(f_0)_{\ast}\del=(f_{n+1})^{\ast}\rho$.

\end{enumerate}
\end{definition}

\begin{definition}\cite[Definition 2.13]{HLN}\label{def1}
 An {\it $n$-exangle} is an object $\Xd$ in $\AE$ that satisfies the listed conditions.
\begin{enumerate}
\item[\rm (1)] The following sequence of functors $\C\op\to{\rm Ab}$ is exact.
$$
\C(-,X_0)\xrightarrow{\C(-,\ d^X_0)}\cdots\xrightarrow{\C(-,\ d^X_n)}\C(-,X_{n+1})\xrightarrow{~\del\ssh~}\E(-,X_0)
$$
\item[\rm (2)] The following sequence of functors $\C\to{\rm Ab}$ is exact.
$$
\C(X_{n+1},-)\xrightarrow{\C(d^X_n,\ -)}\cdots\xrightarrow{\C(d^X_0,\ -)}\C(X_0,-)\xrightarrow{~\del\ush~}\E(X_{n+1},-)
$$
\end{enumerate}
In particular any $n$-exangle is an object in $\AE$.
A {\it morphism of $n$-exangles} simply means a morphism in $\AE$. Thus $n$-exangles form a full subcategory of $\AE$.
\end{definition}

\begin{definition}\cite[Definition 2.22]{HLN}
Let $\s$ be a correspondence which associates a homotopic equivalence class $\s(\del)=[{}_A{X_{\bullet}}_C]$ to each extension $\del={}_A\del_C$. Such $\s$ is called a {\it realization} of $\E$ if it satisfies the following condition for any $\s(\del)=[X_{\bullet}]$ and any $\s(\rho)=[Y_{\bullet}]$.
\begin{itemize}
\item[{\rm (R0)}] For any morphism of extensions $(a,c)\co\del\to\rho$, there exists a morphism $f_{\bullet}\in\CC(X_{\bullet},Y_{\bullet})$ of the form $f_{\bullet}=(a,f_1,\ldots,f_n,c)$. Such $f_{\bullet}$ is called a {\it lift} of $(a,c)$.
\end{itemize}
In such a case, we simple say that \lq\lq$X_{\bullet}$ realizes $\del$" whenever they satisfy $\s(\del)=[X_{\bullet}]$.

Moreover, a realization $\s$ of $\E$ is said to be {\it exact} if it satisfies the following conditions.
\begin{itemize}
\item[{\rm (R1)}] For any $\s(\del)=[X_{\bullet}]$, the pair $\Xd$ is an $n$-exangle.
\item[{\rm (R2)}] For any $A\in\C$, the zero element ${}_A0_0=0\in\E(0,A)$ satisfies
\[ \s({}_A0_0)=[A\ov{\id_A}{\lra}A\to0\to\cdots\to0\to0]. \]
Dually, $\s({}_00_A)=[0\to0\to\cdots\to0\to A\ov{\id_A}{\lra}A]$ holds for any $A\in\C$.
\end{itemize}
Note that the above condition {\rm (R1)} does not depend on representatives of the class $[X_{\bullet}]$.
\end{definition}

\begin{definition}\cite[Definition 2.23]{HLN}
Let $\s$ be an exact realization of $\E$.
\begin{enumerate}
\item[\rm (1)] An $n$-exangle $\Xd$ is called an $\s$-{\it distinguished} $n$-exangle if it satisfies $\s(\del)=[X_{\bullet}]$. We often simply say {\it distinguished $n$-exangle} when $\s$ is clear from the context.
\item[\rm (2)]  An object $X_{\bullet}\in\CC$ is called an {\it $\s$-conflation} or simply a {\it conflation} if it realizes some extension $\del\in\E(X_{n+1},X_0)$.
\item[\rm (3)]  A morphism $f$ in $\C$ is called an {\it $\s$-inflation} or simply an {\it inflation} if it admits some conflation $X_{\bullet}\in\CC$ satisfying $d_0^X=f$.
\item[\rm (4)]  A morphism $g$ in $\C$ is called an {\it $\s$-deflation} or simply a {\it deflation} if it admits some conflation $X_{\bullet}\in\CC$ satisfying $d_n^X=g$.
\end{enumerate}
\end{definition}

\begin{definition}\cite[Definition 2.32]{HLN}
An {\it $n$-exangulated category} is a triplet $(\C,\E,\s)$ of additive category $\C$, additive bifunctor $\E\co\C\op\times\C\to\Ab$, and its exact realization $\s$, satisfying the following conditions.
\begin{itemize}
\item[{\rm (EA1)}] Let $A\ov{f}{\lra}B\ov{g}{\lra}C$ be any sequence of morphisms in $\C$. If both $f$ and $g$ are inflations, then so is $g\circ f$. Dually, if $f$ and $g$ are deflations, then so is $g\circ f$.

\item[{\rm (EA2)}] For $\rho\in\E(D,A)$ and $c\in\C(C,D)$, let ${}_A\langle X_{\bullet},c\uas\rho\rangle_C$ and ${}_A\Yr_D$ be distinguished $n$-exangles. Then $(\id_A,c)$ has a {\it good lift} $f_{\bullet}$, in the sense that its mapping cone gives a distinguished $n$-exangle $\langle M^f_{\bullet},(d^X_0)\sas\rho\rangle$.
 \item[{\rm (EA2$\op$)}] Dual of {\rm (EA2)}.
\end{itemize}
Note that the case $n=1$, a triplet $(\C,\E,\s)$ is a  $1$-exangulated category if and only if it is an extriangulated category, see \cite[Proposition 4.3]{HLN}.
\end{definition}

\begin{example}
From \cite[Proposition 4.34]{HLN} and \cite[Proposition 4.5]{HLN},  we know that $n$-exact categories and $(n+2)$-angulated categories are $n$-exangulated categories.
There are some other examples of $n$-exangulated categories
 which are neither $n$-exact nor $(n+2)$-angulated, see \cite{HLN,LZ,HZZ1}.
\end{example}

The following some lemmas are very useful which are needed later on.

\begin{lemma}\emph{\cite[Proposition 3.6]{HLN}}\label{a2}
\rm Let ${}_A\langle X_{\bullet},\delta\rangle_C$ and ${}_B\langle Y_{\bullet},\rho\rangle_D$ be distinguished $n$-exangles. Suppose that we are given a commutative square
$$\xymatrix{
 X_0 \ar[r]^{{d_0^X}} \ar@{}[dr]|{\circlearrowright} \ar[d]_{a} & X_1 \ar[d]^{b}\\
 Y_0  \ar[r]_{d_0^Y} &Y_1
}
$$
in $\C$. Then the following holds.

{\rm (1)}~ There is a morphism $f_{\bullet}\colon\Xd\to\langle Y_{\bullet},\rho\rangle$ which satisfies $f_0=a$ and $f_1=b$.

{\rm (2)}~ If $X_0=Y_0=A$ and $a=1_A$ for some $A\in\C$, then the above $f_{\bullet}$ can be taken to give a distinguished $n$-exangle $\langle M^f_{\bullet},(d^X_0)\sas\rho\rangle$.

\end{lemma}

Let $(\C,\E,\s)$ be an $n$-exangulated category and  $\F\subseteq\E$ be an additive subfunctor {\rm (see \cite[Definition 3.7]{HLN})}. For a realization $\s$ of $\E$, define
$\s\hspace{-1.3mm}\mid_{\F}$ to be the restriction of $\s$ onto $\F$. Namely, it is defined by $\s\hspace{-1.3mm}\mid_{\F}(\del)=\s(\del)$ for any $\F$-extension $\del$.

\begin{lemma}\label{lem1}{\rm\cite[Proposition 3.16]{HLN}}
Let $(\C,\E,\s)$ be an $n$-exangulated category. For any additive subfunctor $\F\subseteq\E$, the following statements are equivalent.

{\rm (1)}~ $(\C,\F,\s\hspace{-1.2mm}\mid_{\F})$ is an $n$-exangulated category.

{\rm (2)}~ $\s\hspace{-1.2mm}\mid_{\F}$-inflations are closed under composition.

{\rm (3)}~ $\s\hspace{-1.2mm}\mid_{\F}$-deflations are closed under composition.

\end{lemma}

\section{Main result}
In this section, when we say that $\A$ is a subcategory of an additive category $\C$, we always assume that $\C$ is full, and closed under isomorphisms, direct sums and direct summands.

We first recall the notion of $n$-extension closed from \cite{HLN}.
\begin{definition}\cite[Definition 2.34]{HLN}
Let $(\C,\E,\s)$ be an $n$-exangulated category. A subcategory $\A$ of $\C$ is called \emph{$n$-extension closed} if
 for any pair of objects $A$ and $C$ in $\A$ and any $\E$-extension $\delta\in\E(C,A)$, there is a distinguished $n$-exangle $\Xd$ with $X_i$ in $\A$ for $i=1,\cdots,n$.
\end{definition}

Let $(\C,\E,\s)$ be an $n$-exangulated category and $\A$ be an $n$-extension closed subcategory of $\C$. We define $\E_{\A}$ to be the restriction of
$\E$ onto $\A^{\rm op}\times\A$. For any $\delta\in\E_{\A}(C,A)$, take an $\s$-distinguished $n$-exangle $\Xd$ with $X_i$ in $\A$ for $i=1,\cdots,n$. We define $\t(\delta)=[X_{\bullet}]$, where the homotopy  equivalence class is taken in ${\bf{C}}^{n+2}_{({\A;A,C})}$.

Our first main result is the following. The construction gives $n$-exangulated categories which are not
$n$-exact nor $(n+2)$-angulated in general.

\begin{theorem}\label{main1}
Let $(\C,\E,\s)$ be a Krull-Schmidt $n$-exangulated categor and
$\A$ be an $n$-extension closed subcategory of $\C$. Then $(\A,\E_{\A},\t)$
is an $n$-exangulated category.
\end{theorem}
{\bf In order to prove Theorem \ref{main1}, we need some preparations as follows.}

It is obvious that we have the following commutative diagram
$$\xymatrix{
X_{\bullet}\ar[d]^{g_{\bullet}}&0\ar[r]^{}\ar@{}[dr] \ar[d]^{} &0 \ar[r]^{} \ar@{}[dr]\ar[d]^{}&\cdot\cdot\cdot \ar[r]^{} \ar@{}[dr]&0 \ar[r]^{0} \ar@{}[dr]\ar[d]^{}& X \ar[r]^{\id_{X}}\ar@{-->}[ddl]_(0.25){0} \ar@{}[dr]\ar[d]^{}& X \ar[r]^{0} \ar@{-->}[ddl]_(0.25){-\id_{X}}\ar@{}[dr]\ar[d]^{} &0 \ar[r]^{} \ar@{-->}[ddl]_(0.25){0}\ar@{}[dr]\ar[d]^{}&\cdot\cdot\cdot \ar[r]^{} \ar@{}[dr]&0 \ar@{}[dr]\ar[d]
\ar@{-->}[r]^-{0} &\\
0_{\bullet}\ar[d]^{f_{\bullet}}&0\ar[r]^{}\ar@{}[dr] \ar[d]^{} &0\ar[r]^{} \ar@{}[dr]\ar[d]^{}&\cdot\cdot\cdot\ar[r]^{} \ar@{}[dr]&0 \ar[r]^{} \ar@{}[dr]\ar[d]^{}& 0 \ar[r]^{} \ar@{}[dr]\ar[d]^{}& 0 \ar[r]\ar@{}[dr]\ar[d]^{} &0 \ar[r]^{} \ar@{}[dr]\ar[d]^{}&\cdot\cdot\cdot \ar[r]^{} \ar@{}[dr]&{0} \ar@{}[dr]\ar[d]^{} \ar@{-->}[r]^-{0} &\\
X_{\bullet}&0 \ar[r]^{}&0 \ar[r]^{} & \cdot\cdot\cdot\ar[r]^{}  & 0\ar[r]^{0}  & X \ar[r]^{\id_{X}}&X \ar[r]^{0} & 0\ar[r]^{}  &\cdot\cdot\cdot \ar[r]^{}  &{0}\ar@{-->}[r]^-{0}&}
$$
of distinguished $n$-exangles.
It is easy to see that $f_{\bullet}g_{\bullet}\sim\id_{X_{\bullet}}$ and $g_{\bullet}f_{\bullet}\sim\id_{0_{\bullet}}$, hence $$\s({}_00_0)=[0\xrightarrow{}0\xrightarrow{}\cdots\xrightarrow{}0
]=[0\xrightarrow{}0\xrightarrow{}\cdots\xrightarrow{}X\xrightarrow{\id_{X}}X\xrightarrow{}0\xrightarrow{}\cdots\xrightarrow{}0].$$ We denote the split distinguished $n$-exangle $0\xrightarrow{}0\xrightarrow{}\cdots\xrightarrow{}X\xrightarrow{\id_{X}}X\xrightarrow{}0\xrightarrow{}\cdots\xrightarrow{}0\overset{0}{\dashrightarrow}$ in $\C$ by ${\rm{split}}_{i}(X)$, where $\id_{X}$ as the identity morphism in position $i$, and $i\in\{{1,2,\cdots,n-1}\}$.
\medskip

We denote by ${\rm rad}_{\C}$ the Jacobson radical of $\C$. Namely, ${\rm rad}_{\C}$ is an ideal of $\C$ such that ${\rm rad}_{\C}(A, A)$
coincides with the Jacobson radical of the endomorphism ring ${\rm End}(A)$ for any $A\in\C$.
\begin{definition} When $n\geq2$, a distinguished $n$-exangle in $\C$ of the form
$$A_{\bullet}:~~A_0\xrightarrow{\alpha_0}A_1\xrightarrow{\alpha_1}A_2\xrightarrow{\alpha_2}\cdots\xrightarrow{\alpha_{n-2}}A_{n-1}
\xrightarrow{\alpha_{n-1}}A_n\xrightarrow{\alpha_n}A_{n+1}\overset{}{\dashrightarrow}$$
is minimal if $\alpha_1,\alpha_2,\cdots,\alpha_{n-1}$ are in $\rad_{\C}$.

\end{definition}
The following lemma shows that $\E$-extension in an equivalence class can be chosen in a minimal way in a Krull-Schmidt $n$-exangulated category.

\begin{lemma}\rm\label{ml} Let $\C$ be a Krull-Schmidt $n$-exangulated category, $A_0,A_{n+1}\in\C$. Then for every equivalence class of $\E$-extension of $A_0$ by $A_{n+1}$, there exists a representation
$$A_{\bullet}:~~A_0\xrightarrow{\alpha_0}A_1\xrightarrow{\alpha_1}A_2\xrightarrow{\alpha_2}\cdots\xrightarrow{\alpha_{n-2}}A_{n-1}
\xrightarrow{\alpha_{n-1}}A_n\xrightarrow{\alpha_n}A_{n+1}\overset{}{\dashrightarrow}$$
such that $\alpha_1,\alpha_2,\cdots,\alpha_{n-1}$ are in $\rad_{\C}$. Moreover, $A_{\bullet}$ is a direct summand of every other equivalent $\E$-extension.
\end{lemma}
\proof
Let $\delta\in\E(A_{n+1},A_0)$, there exists a distinguished $n$-exangle in $\C$ of the form $$X_{\bullet}:~~A_0\xrightarrow{\beta_0}X_1\xrightarrow{\beta_1}X_2\xrightarrow{\beta_2}\cdots\xrightarrow{\beta_{n-2}}X_{n-1}
\xrightarrow{\beta_{n-1}}X_n\xrightarrow{\beta_n}A_{n+1}\overset{\delta}{\dashrightarrow}.$$
Suppose $\beta_i$ is not in the  $\rad_{\C}$ for some $i\in\{{1,2,\cdots,n-1}\}$. By definition there is a nonzero object $Z$ and $g:Z\rightarrow X_i$ and $h:X_{i+1}\rightarrow Z$ such that $h\beta_ig$ is an isomorphism. Then we have the following commutative diagram
$$\xymatrix{
0\ar[r]^{}\ar@{}[dr] \ar[d]^{} &0 \ar[r]^{} \ar@{}[dr]\ar[d]^{}&\cdot\cdot\cdot \ar[r]^{} \ar@{}[dr]&0 \ar[r]^{} \ar@{}[dr]\ar[d]^{}& Z\ar@{=}[r]^{} \ar@{}[dr]\ar[d]^{g}& Z \ar[r]^{} \ar@{}[dr]\ar[d]^{\beta_ig} &0 \ar[r]^{} \ar@{}[dr]\ar[d]^{}&\cdot\cdot\cdot \ar[r]^{} \ar@{}[dr]&0 \ar@{}[dr]\ar[d]^{} \ar@{-->}[r]^-{0} &\\A_0\ar[r]^{\beta_0}\ar@{}[dr] \ar[d]^{} &X_1\ar[r]^{\beta_1} \ar@{}[dr]\ar[d]^{}&\cdot\cdot\cdot\ar[r]^{} \ar@{}[dr]&X_{i-1} \ar[r]^{\beta_{i-1}} \ar@{}[dr]\ar[d]^{}& X_{i} \ar[r]^{\beta_{i}} \ar@{}[dr]\ar[d]|{(h\beta_ig)^{-1}h\beta_i}& X_{i+1}\ar[r]^{\beta_{i+1}}\ar@{}[dr]\ar[d]|{(h\beta_ig)^{-1}h} &X_{i+2}\ar[r]^{\beta_{i+2}} \ar@{}[dr]\ar[d]^{}&\cdot\cdot\cdot \ar[r]^{\beta_{n}} \ar@{}[dr]&{A_{n+1}} \ar@{}[dr]\ar[d]^{} \ar@{-->}[r]^-{} &\\
0 \ar[r]^{}&0 \ar[r]^{} & \cdot\cdot\cdot\ar[r]^{}  & 0\ar[r]^{}  & Z \ar@{=}[r]^{}&Z \ar[r]^{} & 0\ar[r]^{}  &\cdot\cdot\cdot \ar[r]^{}  &{0} \ar@{-->}[r]^-{0} &}
$$
of distinguished $n$-exangles. Hence the split distinguished $n$-exangle ${\rm{split}}_{i}(Z)$ is a direct summand of $X_{\bullet}$. Let $X_{\bullet}={\rm{split}}_{i}(Z)\oplus Y_{\bullet}$ as distinguished $n$-exangles, where $Y_{\bullet}$ is a complement of $X_{\bullet}$. If $Y_{\bullet}$ is not a minimal distinguished $n$-exangle, we can repeat this procedure. Since $\C$ is a Krull-Schmidt category, the $X_{\bullet}$ has the descending chain condition on direct summand. This means that after finitely many steps this procedure must stop. Notice that ${\rm{split}}_{i}(Z)\oplus{\rm{split}}_{i}(Z^{\prime})={\rm{split}}_{i}(Z\oplus Z^{\prime})$ for $i\in\{{1,2,\cdots,n-1}\}$ and $Z, Z^{\prime}\in\C$ by Proposition 3.3 in \cite{HLN}. Hence $X_{\bullet}$ is the direct sum of a minimal distinguished $n$-exangle and split distinguished $n$-exangles. In particular, this also shows a minimal distinguished $n$-exangle is exists. Finally, we only need to show that the minimal distinguished $n$-exangle $A_{\bullet}$ and the distinguished $n$-exangle $X_{\bullet}$ are homotopy equivalence. By $\rm(R0)$ we have the following commutative diagram of distinguished $n$-exangles

$$\xymatrix{
A_{\bullet}\ar[d]^{f_{\bullet}}:& A_0\ar[r]^{\alpha_0}\ar@{}[dr] \ar@{=}[d]^{} &A_1 \ar[r]^{\alpha_1} \ar@{}[dr]\ar@{-->}[d]^{f_1}&A_2\ar[r]^{\alpha_2} \ar@{}[dr]\ar@{-->}[d]^{f_2}&\cdot\cdot\cdot  \ar[r]^{\alpha_{n-1}} \ar@{}[dr]& A_n\ar[r]^{\alpha_{n}} \ar@{}[dr]\ar@{-->}[d]^{f_n}& A_{n+1}  \ar@{}[dr] \ar@{=}[d]^{}\ar@{-->}[r]^-{} &\\X_{\bullet}\ar[d]^{g_{\bullet}}:& A_0\ar[r]^{\beta_0}\ar@{}[dr]\ar@{=}[d]^{} &X_1\ar[r]^{\beta_1} \ar@{}[dr]\ar@{-->}[d]^{g_1}&X_2\ar[r]^{\beta_2} \ar@{}[dr]\ar@{-->}[d]^{g_2}&\cdot\cdot\cdot \ar[r]^{\beta_{n-1}} \ar@{}[dr]& X_n \ar[r]^{\beta_{n}}\ar@{-->}[d]^{g_n} \ar@{}[dr]& X_{n+1} \ar@{}[dr]\ar@{=}[d]^{} \ar@{-->}[r]^-{} &\\
A_{\bullet}& A_0\ar[r]^{\alpha_0}&A_1 \ar[r]^{\alpha_1} & A_2\ar[r]^{\alpha_2}  & \cdot\cdot\cdot\ar[r]^{\alpha_{n-1}}  & A_n \ar[r]^{\alpha_{n}}&A_{n+1} \ar@{-->}[r]^-{} &}$$
Since $\alpha_{i+1}$ is a weak cokernel of $\alpha_{i}$ and $\beta_{i+1}$ is a weak cokernel of $\beta_{i}$ by Definition \ref{def1}, where $i=0,1,\cdots,n-1$. Notice that $f_{0}g_{0}=\id_{A_0}$ and $g_{0}f_{0}=\id_{A_0}$. Then we have $f_{\bullet}g_{\bullet}\sim\id_{X_{\bullet}}$ and $g_{\bullet}f_{\bullet}\sim\id_{A_{\bullet}}$ by the Comparison Lemma 2.1 in \cite{Ja}.

\qed

\begin{remark} Let $\C$ be a Krull-Schmidt $n$-exangulated category. By the Krull-Schmidt property of $\C$, every minimal distinguished $n$-exangle in each equivalence class is unique up to isomorphism.
\end{remark}

\begin{construction}\label{con}
For any commutative diagram of distinguished $n$-exangles
\begin{equation}\label{t3}
\begin{array}{l}
\xymatrix{
X_{\bullet}:& X_0\ar[r]^{f_0}\ar@{}[dr] |{\circlearrowright}\ar@{=}[d]^{} &Y_0 \ar[r]^{f_1} \ar@{}[dr]\ar[d]^{g_0}&X_2\ar[r]^{f_2} \ar@{}[dr]&\cdot\cdot\cdot  \ar[r]^{f_{n-1}} \ar@{}[dr]& X_n\ar[r]^{f_{n}} \ar@{}[dr]& X_{n+1}  \ar@{}[dr] \ar@{-->}[r]^-{\delta} &\\Z_{\bullet}:& X_0\ar[r]^{h_0}\ar@{}[dr]|{\circlearrowright} \ar[d]^{f_0} &Z_0\ar[r]^{h_1} \ar@{}[dr]\ar@{=}[d]^{}&Z_2\ar[r]^{h_2} \ar@{}[dr]&\cdot\cdot\cdot \ar[r]^{h_{n-1}} \ar@{}[dr]& Z_n \ar[r]^{h_{n}} \ar@{}[dr]& Z_{n+1} \ar@{}[dr] \ar@{-->}[r]^-{\eta} &\\
Y_{\bullet}:& Y_0\ar[r]^{g_0}&Z_0 \ar[r]^{g_1} & Y_2\ar[r]^{g_2}  & \cdot\cdot\cdot\ar[r]^{g_{n-1}}  & Y_n \ar[r]^{g_{n}}&Y_{n+1} \ar@{-->}[r]^-{\xi} &}
\end{array}
\end{equation}
by Lemma \ref{a2}, there exists a morphism $\varphi_{\bullet}\colon\Xd\to\langle Z_{\bullet},\eta\rangle$ which satisfies $\varphi_0=\id_{X_0}$ and $\varphi_1=g_0$. That is,
there exists the following commutative diagram of distinguished $n$-exangles
$$\xymatrix{
X_0\ar[r]^{f_0}\ar@{}[dr] \ar@{=}[d]^{} &Y_0 \ar[r]^{f_1} \ar@{}[dr]\ar[d]^{g_0}&X_2\ar[r]^{f_2} \ar@{-->}[d]^{\varphi_2}\ar@{}[dr]&\cdot\cdot\cdot  \ar[r]^{f_{n-1}} \ar@{}[dr]& X_n\ar[r]^{f_{n}} \ar@{}[dr]\ar@{-->}[d]^{\varphi_n}& X_{n+1}  \ar@{}[dr] \ar@{-->}[d]^{\varphi_{n+1}}\ar@{-->}[r]^-{\delta} &\\
X_0\ar[r]^{h_0}&Z_0 \ar[r]^{h_1} & Z_2\ar[r]^{h_2}  & \cdot\cdot\cdot\ar[r]^{h_{n-1}}  & Z_n \ar[r]^{h_{n}}&Z_{n+1} \ar@{-->}[r]^-{\eta} &}$$
such that
$$M_{\bullet}:Y_0\xrightarrow{\left[
                                            \begin{smallmatrix}
                                              -f_1 \\
                                              g_0 \\
                                            \end{smallmatrix}
                                          \right]}
 X_2\oplus Z_0\xrightarrow{\left[
                                            \begin{smallmatrix}
                                              -f_{2} & 0 \\
                                              \varphi_{2} & h_1
                                            \end{smallmatrix}
                                          \right]}
X_3\oplus Z_2\xrightarrow{\left[
                                            \begin{smallmatrix}
                                              -f_{3} & 0 \\
                                             \varphi_{3} & h_2
                                            \end{smallmatrix}
                                          \right]}
\cdots\xrightarrow{\left[
                                            \begin{smallmatrix}
                                              -f_{n} & 0 \\
                                              \varphi_{n} & h_{n-1} \\
                                            \end{smallmatrix}
                                          \right]\ }
X_{n+1}\oplus Z_{n}\xrightarrow{\left[\begin{smallmatrix}
                                            \varphi_{n+1}, & h_{n} \\
                                            \end{smallmatrix}
                                          \right]}
Z_{n+1}\overset{(f_0){_{\ast}}\eta}{\dashrightarrow}$$
is a distinguished $n$-exangle.

Similarly, we have the following commutative diagram of distinguished $n$-exangles
  $$\xymatrix@C=1.3cm@R=1.2cm{
Y_0\ar[r]^{\left[
                                            \begin{smallmatrix}
                                              -f_1 \\
                                              g_0 \\
                                            \end{smallmatrix}
                                          \right]}\ar@{}[dr] \ar@{=}[d]^{} &X_2\oplus Z_0 \ar[r]^{\left[
                                            \begin{smallmatrix}
                                              -f_{2} & 0 \\
                                              \varphi_{2} & h_1
                                            \end{smallmatrix}
                                          \right]} \ar@{}[dr]\ar[d]^{\left[\begin{smallmatrix}
                                           0 & 1 \\
                                            \end{smallmatrix}
                                          \right]}&X_3\oplus Z_2\ar[r]^{\left[
                                            \begin{smallmatrix}
                                              -f_{3} & 0 \\
                                             \varphi_{3} & h_2
                                            \end{smallmatrix}
                                          \right]} \ar@{-->}[d]^{\psi_2}\ar@{}[dr]&\cdot\cdot\cdot  \ar[r]^{\left[
                                            \begin{smallmatrix}
                                              -f_{n} & 0 \\
                                              \varphi_{n} & h_{n-1} \\
                                            \end{smallmatrix}
                                          \right]} \ar@{}[dr]& X_{n+1}\oplus Z_{n}\ar[r]^{\left[\begin{smallmatrix}
                                            \varphi_{n+1}, & h_{n} \\
                                            \end{smallmatrix}
                                          \right]} \ar@{}[dr]\ar@{-->}[d]^{\psi_n}& Z_{n+1} \ar@{}[dr] \ar@{-->}[d]^{\psi_{n+1}}\ar@{-->}[r]^-{(f_0){_{\ast}}\eta}  &\\
Y_0\ar[r]^{g_0}&Z_0 \ar[r]^{g_1} & Y_2\ar[r]^{g_2}  & \cdot\cdot\cdot\ar[r]^{g_{n-1}}  & Y_n \ar[r]^{g_{n}}&Y_{n+1} \ar@{-->}[r]^-{\xi}  &}
$$
where $\psi_{i}=\left[\begin{smallmatrix}
                                            m_{i} & k_{i} \\
                                            \end{smallmatrix}
                                          \right], i=2,\cdot\cdot\cdot,n$, such that

$$D_{\bullet}:X_2\oplus Z_0\xrightarrow{\left[
                                                                            \begin{smallmatrix}
                                                                              f_2 & 0 \\
                                                                              -\varphi_2 & -h_1 \\
                                                                             0 & 1\\
                                                                            \end{smallmatrix}
                                                                          \right]} X_3\oplus Z_2\oplus Z_0\xrightarrow{\left[
                                                                            \begin{smallmatrix}
                                                                              f_3 & 0 & 0 \\
                                                                              -\varphi_3 & -h_2 & 0 \\
                                                                              m_2 & k_2 & g_1\\
                                                                            \end{smallmatrix}
                                                                          \right]} X_4\oplus Z_3\oplus Y_2\xrightarrow{\left[
                                                                            \begin{smallmatrix}
                                                                               f_4 & 0 & 0 \\
                                                                              -\varphi_4 & -h_3 & 0 \\
                                                                              m_3 & k_3 & g_2\\
                                                                            \end{smallmatrix}
                                                                          \right]} \cdots\xrightarrow{\left[
                                                                            \begin{smallmatrix}
                                                                              f_n & 0 & 0 \\
                                                                              -\varphi_n & -h_{n-1} & 0 \\
                                                                              m_{n-1} & k_{n-1} & g_{n-2}\\
                                                                            \end{smallmatrix}
                                                                          \right]}
$$
 $$ X_{n+1}\oplus Z_{n}\oplus Y_{n-1}\xrightarrow{\left[
                                                                            \begin{smallmatrix}
                                                                              -\varphi_{n+1} & -h_{n} & 0 \\
                                                                              m_{n} & k_{n} & g_{n-1}\\
                                                                            \end{smallmatrix}
                                                                          \right]} Z_{n+1}\oplus Y_{n}\xrightarrow{[\psi_{n+1}, \ g_{n}]} Y_{n+1}\overset{{\left[
                                            \begin{smallmatrix}
                                              -f_1 \\
                                              g_0 \\
                                            \end{smallmatrix}
                                          \right]}_{\ast}{\xi}}{\dashrightarrow}$$
is a distinguished $n$-exangle.
\end{construction}

The following lemma holds in any additive Krull-Schmidt category $\C$.

\begin{lemma}\label{matrix}{\rm\cite[Lemma 2.7]{K}}
Suppose $f\colon X_{1}\oplus X_{2}\rightarrow Y_{1}\oplus Y_{2}$ is a morphism in any additive Krull-Schmidt category $\C$ with $\pi_{1}f\iota_{1}\in\rad_{\C}$, where $\iota_{1}:X_{1}\rightarrow X_{1}\oplus X_{2}$ is the canonical inclusion and $\pi_{1}:Y_{1}\oplus Y_{2}\rightarrow Y_{1}$ is the canonical projection. Then all objects $Z$ which have $g:Z\rightarrow X_{1}\oplus X_{2}$ and $h:Y_{1}\oplus Y_{2}\rightarrow Z$ such that  $hfg\in \Aut~ Z$ satisfy $Z\in \add(X_{2}\oplus Y_{2})$.
\end{lemma}

The following lemma  was
proved in \cite[Lemma 2.4]{K}, when $\C$ is an $(n+2)$-angulated category. However, it can be
extended to our setting.

\begin{lemma}\label{Key}
Suppose we have a commutative diagram as in {\rm (\ref{t3})} for  $n\geq3$ with the second row ${}_{X_{0}}\langle Z_{\bullet},\eta\rangle_{Z_{n+1}}$  a minimal distinguished $n$-exangle.
For ${\left[
                                            \begin{smallmatrix}
                                              -f_1 \\
                                              g_0 \\
                                            \end{smallmatrix}
                                          \right]}_{\ast}{\xi}\in\E( Y_{n+1},X_2\oplus Z_0)$, let

$$W_{\bullet}:~~X_2\oplus Z_0\xrightarrow{}W_1\xrightarrow{}W_2\xrightarrow{}\cdots\xrightarrow{}W_{n-1}
\xrightarrow{}W_n\xrightarrow{} Y_{n+1}\overset{{\left[
                                            \begin{smallmatrix}
                                              -f_1 \\
                                              g_0 \\
                                            \end{smallmatrix}
                                          \right]}_{\ast}{\xi}}{\dashrightarrow}$$
be a minimal distinguished $n$-exangle. Then there exist objects ${\{N_k\}^{n-1}_{k=1}}$ satisfying

{\rm (1)}~$N_1\in \add(X_{3}\oplus Z_{0}\oplus X_{4}\oplus Y_{2})$,

{\rm (2)}~$N_k\in \add(X_{k+2}\oplus Y_{k}\oplus X_{k+3}\oplus Y_{k+1})$~ for $k=2,3,\cdots,n-2$,

{\rm (3)}~$N_{n-1}\in \add(X_{n+1}\oplus  Y_{n-1}\oplus Y_{n})$.\\
Moreover, the following relations hold

{\rm (4)}~$X_{3}\oplus Z_{2}\oplus Z_{0}=W_{1}\oplus N_{1}$,

{\rm (5)}~$X_{k+2}\oplus Z_{k+1}\oplus Y_{k}=N_{k-1}\oplus W_{k}\oplus N_{k}$~ for $k=2,3,\cdots,n-2$,

{\rm (6)}~$Z_{n+1}\oplus Y_{n}=N_{n-1}\oplus W_{n}$,\\
which reveals the difference between the distinguished $n$-exangle ${}_{X_{2}\oplus Z_{0}}\langle D_{\bullet},{\left[
                                            \begin{smallmatrix}
                                              -f_1 \\
                                              g_0 \\
                                            \end{smallmatrix}
                                          \right]}_{\ast}{\xi}\rangle_{Y_{n+1}}$ and the minimal distinguished $n$-exangle ${}_{X_{2}\oplus Z_{0}}\langle W_{\bullet},{\left[
                                            \begin{smallmatrix}
                                              -f_1 \\
                                              g_0 \\
                                            \end{smallmatrix}
                                          \right]}_{\ast}{\xi}\rangle_{Y_{n+1}}$.

\end{lemma}

\proof
By Lemma \ref{ml}, we have that the distinguished $n$-exangle $D_{\bullet}$ is the direct sum of the minimal distinguished $n$-exangle $W_{\bullet}$ and split distinguished $n$-exangle $\rm{split}_{1}(N_{1}),\cdots,\rm{split}_{n-1}(N_{n-1})$. This choice of $ N_{1},\cdots,N_{n-1}$ establishes $(4)-(6)$ of Lemma.

It remains to show that $(1)-(3)$ of Lemma hold for this choice of $ N_{1},\cdots,N_{n-1}$. We only show $(1)$, the other statements are analogous. Let $\zeta=\left[
                                                                            \begin{smallmatrix}
                                                                              f_3 & 0 & 0 \\
                                                                              -\varphi_3 & -h_2 & 0 \\
                                                                              m_2 & k_2 & g_1\\
                                                                            \end{smallmatrix}
                                                                          \right]$,
since the split distinguished $n$-exangle ${\rm{split}}_{1}(N_{1})$ is a direct summand of the distinguished $n$-exangle $D_{\bullet}$, we have that $\id_{N_{1}}$ factors through $\zeta$. Since $-h_2 $ is in $\rad_{\C}$ by assumption, we have that $N_1\in \add(X_{3}\oplus Z_{0})\oplus( X_{4}\oplus Y_{2})$ by Lemma \ref{matrix}.
\qed

\begin{lemma}\label{com}
If both $f_0$ and $g_0$ are $\t$-inflations, then so is $h_0=g_0f_0$.

\end{lemma}
\proof Let $f_0:X_{0}\rightarrow Y_{0}$ and $g_0:Y_{0}\rightarrow Z_{0}$ be $\t$-inflations and $h_0=g_0f_0:X_{0}\rightarrow  Z_{0}$. Since $f_0$ and $g_0$ are $\t$-inflations, we have the following two $\s$-distinguished $n$-exangles
$$X_{\bullet}:~~X_0\xrightarrow{f_0}Y_0\xrightarrow{f_1}X_2\xrightarrow{f_2}\cdots\xrightarrow{f_{n-2}}X_{n-1}
\xrightarrow{f_{n-1}}X_n\xrightarrow{f_n}X_{n+1}\overset{\delta}{\dashrightarrow}$$ and
$$Y_{\bullet}:~~Y_0\xrightarrow{g_0}Z_0\xrightarrow{g_1}Y_2\xrightarrow{g_2}\cdots\xrightarrow{g_{n-2}}Y_{n-1}
\xrightarrow{g_{n-1}}Y_n\xrightarrow{g_n}Y_{n+1}\overset{\xi}{\dashrightarrow}, $$
where $X_2,\cdots,X_{n+1};Y_2,\cdots,Y_{n+1}\in\A$. By {\rm (EA1)}, we know that $h_0$ is an $\s$-inflation, we can choice a minimal $\s$-distinguished $n$-exangle by Lemma \ref{ml}
$$Z_{\bullet}:~~X_0\xrightarrow{h_0}Z_0\xrightarrow{h_1}Z_2\xrightarrow{h_2}\cdots\xrightarrow{h_{n-2}}Z_{n-1}
\xrightarrow{h_{n-1}}Z_n\xrightarrow{h_n}Z_{n+1}\overset{\eta}{\dashrightarrow}, $$ where $Z_2,\cdots,Z_{n+1}\in\C$. We only need to show that $Z_2,\cdots,Z_{n+1}\in\A$. We use the same notation as Construction \ref{con} and Lemma \ref{Key}. For ${\left[
                                            \begin{smallmatrix}
                                              -f_1 \\
                                              g_0 \\
                                            \end{smallmatrix}
                                          \right]}_{\ast}{\xi}\in\E( Y_{n+1},X_2\oplus Z_0)$ with $ Y_{n+1},X_2\oplus Z_0\in\A$, since $\A$ is $n$-extension closed, there is a $\s$-distinguished $n$-exangle
$$S_{\bullet}:~X_2\oplus Z_0\xrightarrow{}S_1\xrightarrow{}S_2\xrightarrow{}\cdots\xrightarrow{}S_{n-1}
\xrightarrow{}S_n\xrightarrow{}Y_{n+1}\overset{{\left[
                                            \begin{smallmatrix}
                                              -f_1 \\
                                              g_0 \\
                                            \end{smallmatrix}
                                          \right]}_{\ast}{\xi}}{\dashrightarrow}, $$
where $S_1,\cdots,S_n\in\A$. Since the minimal $\s$-distinguished $n$-exangle $W_{\bullet}$ is a direct summand of the $\s$-distinguished $n$-exangle $S_{\bullet}$, then $W_1,\cdots,W_n\in\A$. Notice that $N_1,\cdots,N_{n-1}\in\A$ by Lemma \ref{Key} (1)-(3). Therefore, $Z_2,\cdots,Z_{n+1}\in\A$ by Lemma \ref{Key} (4)-(6). This shows that $h_0$ is a $\t$-inflation.
\qed
\medskip

{\bf Now we are ready to prove Theorem \ref{main1}.}
\proof It is straightforward to verify that $\E_{\A}$ is an additive subfunctor of $\E$.
By Lemma \ref{com}, we know that $\t$-inflations are closed under composition.
By Lemma \ref{lem1}, we have that $(\A,\E_{\A},\t)$
is an $n$-exangulated category. \qed
\medskip

By applying Theorem \ref{main1} to $(n+2)$-angulated category, we have the following.
\begin{corollary}\rm\cite[Theorem 3.4]{Z}
Let $(\C,\Sigma,\Theta)$ be a Krull-Schmidt $(n+2)$-angulated category and
$\A$ be an $n$-extension closed subcategory of $\C$. Then $(\A,\E_{\A},\t)$
is an $n$-exangulated category.
\end{corollary}

By applying Theorem \ref{main1} to $n$-exact category, we have the following.

\begin{corollary}
Let $\B$ be a Krull-Schmidt $n$-exact category and
$\A$ be an $n$-extension closed subcategory of $\B$. Then $(\A,\E_{\A},\t)$
is an $n$-exangulated category.
\end{corollary}

In Theorem \ref{main1}, we know that $(\A,\E_{\A},\t)$ is an $n$-exangulated category.
Next we  give a sufficient
condition on when an $n$-exangulated category $\A$ is an $n$-exact category.

\begin{definition}\label{def2}\cite[Definition 3.2]{LZ}
Let $(\C,\E,\s)$ be an $n$-exangulated category.
\begin{itemize}
\item[(1)] An object $P\in\C$ is called \emph{projective} if for any distinguished $n$-exangle
$$A_0\xrightarrow{\alpha_0}A_1\xrightarrow{\alpha_1}A_2\xrightarrow{\alpha_2}\cdots\xrightarrow{\alpha_{n-2}}A_{n-1}
\xrightarrow{\alpha_{n-1}}A_n\xrightarrow{\alpha_n}A_{n+1}\overset{\delta}{\dashrightarrow}$$
and any morphism $c$ in $\C(P,A_{n+1})$, there exists a morphism $b\in\C(P,A_n)$ satisfying $\alpha_n\circ b=c$.
We denote the full subcategory of projective objects in $\C$ by $\P$.
Dually, the full subcategory of injective objects in $\C$ is denoted by $\I$.

\item[(2)] We say that $\C$ {\it has enough  projectives} if
for any object $C\in\C$, there exists a distinguished $n$-exangle
$$B\xrightarrow{\alpha_0}P_1\xrightarrow{\alpha_1}P_2\xrightarrow{\alpha_2}\cdots\xrightarrow{\alpha_{n-2}}P_{n-1}
\xrightarrow{\alpha_{n-1}}P_n\xrightarrow{\alpha_n}C\overset{\delta}{\dashrightarrow}$$
satisfying $P_1,P_2,\cdots,P_n\in\P$. We can define the notion of having \emph{enough injectives} dually.
\end{itemize}
\end{definition}

\begin{definition}
Let $\C^{\prime}$ and $\C^{\prime\prime}$ be two subcategories of an $n$-exangulated category $\C$.
\begin{itemize}
\item[(a)] Denote by $\CoCone(\C^{\prime},\C^{\prime\prime})$ the subcategory

$\{A\in\C\mid$ there exists a distinguished $n$-exangle $A\xrightarrow{}C^{\prime}_1\xrightarrow{}C^{\prime}_2\xrightarrow{}\cdots\xrightarrow{}C^{\prime}_{n-1}
\xrightarrow{}C^{\prime}_n\xrightarrow{}$

~~~~~~~~~~~~$C^{\prime\prime}\overset{}{\dashrightarrow}$ with $C^{\prime}_1,\cdots,C^{\prime}_n\in\C^{\prime}$ and $C^{\prime\prime}\in\C^{\prime\prime}\}$.

\item[(b)] Denote by $\Cone(\C^{\prime},\C^{\prime\prime})$ the subcategory

$\{A\in\C\mid$ there exists a distinguished $n$-exangle $C^{\prime}\xrightarrow{}C^{\prime\prime}_1\xrightarrow{}C^{\prime\prime}_2\xrightarrow{}\cdots\xrightarrow{}C^{\prime\prime}_{n-1}
\xrightarrow{}C^{\prime\prime}_n\xrightarrow{}$

~~~~~~~~~~~~$A\overset{}{\dashrightarrow}$ with $C^{\prime\prime}_1,\cdots,C^{\prime\prime}_n\in\C^{\prime\prime}$ and $C^{\prime}\in\C^{\prime}\}$.

\item[(c)] Define $\Omega\C^{\prime}=\CoCone(\P,\C^{\prime})$.
We write an object $B$ in the form $\Omega C$ if it admits a distinguished $n$-exangle
$$B\xrightarrow{\alpha_0}P_1\xrightarrow{\alpha_1}P_2\xrightarrow{\alpha_2}\cdots\xrightarrow{\alpha_{n-2}}P_{n-1}
\xrightarrow{\alpha_{n-1}}P_n\xrightarrow{\alpha_n}C\overset{\delta}{\dashrightarrow},$$
where $P_1,P_2,\cdots,P_n\in\P$.

\item[(d)] Define $\Sigma\C^{\prime}=\Cone(\C^{\prime},\I)$.
We write an object $C^{\prime}$ in the form $\Sigma B^{\prime}$ if it admits a distinguished $n$-exangle
$$B^{\prime}\xrightarrow{\beta_0}I_1\xrightarrow{\beta_1}I_2\xrightarrow{\beta_2}\cdots\xrightarrow{\beta_{n-2}}I_{n-1}
\xrightarrow{\beta_{n-1}}I_n\xrightarrow{\beta_n}C^{\prime}\overset{\eta}{\dashrightarrow}$$
where $I_1,I_2,\cdots,I_n\in\I$.
\end{itemize}
\end{definition}

\begin{lemma}\label{Key2}
Let $(\C,\E,\s)$ be an $n$-exangulated category and
$$A_{\bullet}:~A_0\xrightarrow{{\left[
                                            \begin{smallmatrix}
                                              f_0 \\
                                              g_0 \\
                                            \end{smallmatrix}
                                          \right]}}A_1\oplus B_1\xrightarrow{[f_1,\hspace{0.8mm}g_{1}] }A_2\xrightarrow{f_2}A_3\xrightarrow{f_3}\cdots\xrightarrow{f_{n-1}}A_n\xrightarrow{f_n}A_{n+1}\overset{}{\dashrightarrow}.$$
is a distinguished $n$-exangle in $\C$.
If $g_0=0$, then $A_{\bullet}\simeq A'_{\bullet}\oplus B_{\bullet}$, where
$$A'_{\bullet}:~A_0\xrightarrow{f_0}A_1\xrightarrow{f_{11}}A'_2\xrightarrow{f_{21}}A_3\xrightarrow{f_3}\cdots\xrightarrow{f_{n-1}}A_n\xrightarrow{f_n}A_{n+1}$$
and
$$B_{\bullet}:~0\xrightarrow{}B_1\xrightarrow{~1~}B_1\xrightarrow{}0\xrightarrow{}\cdots\xrightarrow{}0\xrightarrow{}0\xrightarrow{}0$$
is two distinguished $n$-exangles.
\end{lemma}
\proof The proof is an adaption of \cite[Lemma 2.5]{L}.
We first show that $g_0=0$ implies that $g_1$ is a section. In fact, since $B_{\bullet}$ is a split distinguished $n$-exangle,  if $g_0=0$, then we have the following commutative diagram of distinguished $n$-exangles by Lemma \ref{a2}
$$\xymatrix{
A_0\ar[r]^{{\left[
                                            \begin{smallmatrix}
                                              f_0 \\
                                              g_0 \\
                                            \end{smallmatrix}
                                          \right]}\quad\;\;}\ar[d] & A_1\oplus B_1\ar[r]^{\quad[f_1,\hspace{0.8mm}g_{1}]}\ar[d]^{[0,\hspace{0.8mm}{1}]} &  A_2\ar[r]^{f_2}\ar@{-->}[d]^{g^{\prime}_{2}} & \cdots \ar[r]^{f_{n-1}} & A_n\ar[r]^{f_n}\ar@{-->}[d] & A_{n+1}\ar@{-->}[d] \ar@{-->}[r]^-{} &\\
0\ar[r] &B_1\ar[r]^{1} & B_1\ar[r] & \cdots \ar[r] & 0 \ar[r] & 0\ar@{-->}[r]^-{} &\\
}$$
 Thus there exists a morphism $g^{\prime}_{2}:A_2\rightarrow B_1$ such that $g^{\prime}_{2}f_1=0$ and $g^{\prime}_{2}g_1=1$. Therefore, $g_1$ is a section.

Now we can assume that
$$A_\bullet=(A_0\xrightarrow{\left[
                         \begin{smallmatrix}
                           f_0 \\
                           0 \\
                         \end{smallmatrix}
                       \right]}A_1\oplus B_1\xrightarrow{\left[
                                                           \begin{smallmatrix}
                                                             f_{11} & 0 \\
                                                             f_{12} & g_{11} \\
                                                           \end{smallmatrix}
                                                         \right]
                       }A^{\prime}_2\oplus g_1(B_1)\xrightarrow{[f_{21},\hspace{0.8mm}f_{22}]} A_3\xrightarrow{f_3}\cdots\xrightarrow{f_{n-1}} A_n\xrightarrow{f_n}A_{n+1}\overset{}{\dashrightarrow})$$
where $g_{11}$ is an isomorphism and $f_{22}=0$. Since $[f_{12},\hspace{0.8mm}0]\left[
                         \begin{smallmatrix}
                           f_0 \\
                           0 \\
                         \end{smallmatrix}
                       \right]=f_{12}f_0=0$, there exists a morphism $[a,\hspace{0.8mm}b]: A^{\prime}_2\oplus g_1(B_1)\rightarrow g_1(B_1)$, such that $[f_{12},\hspace{0.8mm}0]=[a,\hspace{0.8mm}b]\left[
                                                           \begin{smallmatrix}
                                                             f_{11} & 0 \\
                                                             f_{12} & g_{11} \\
                                                           \end{smallmatrix}
                                                         \right]$.
Thus $b=0$ and $f_{12}=af_{11}$. The following commutative diagram
$$\xymatrixcolsep{4pc}\xymatrixrowsep{3pc}\xymatrix{A_1\oplus B_1\ar[r]^{\left[
                                                           \begin{smallmatrix}
                                                           f_{11} & 0 \\
                                                              0 & 1 \\
                                                           \end{smallmatrix}
                                                         \right]
                       }\ar@{=}[d] & A^{\prime}_2 \oplus B_1\ar[r]^{[f_{21},\hspace{0.8mm}0]}\ar[d]^{\left[
                                                                                          \begin{smallmatrix}
                                                                                            1 & 0 \\
                                                                                             a & g_{11} \\
                                                                                          \end{smallmatrix}
                                                                                        \right]
                       } & A_3\ar@{=}[d]\\
A_1\oplus B_1\ar[r]^{\left[
                                                           \begin{smallmatrix}
                                                             f_{11} & 0 \\
                                                            f_{12} & g_{11} \\
                                                           \end{smallmatrix}
                                                         \right]
                       }& A^{\prime}_2\oplus g_1(B_1)\ar[r]^{[f_{21},\hspace{0.8mm}0]} & A_3\\}$$
implies that $ A'_\bullet\oplus B_\bullet\cong A_\bullet$. Thus we have $A'_\bullet$ is a distinguished $n$-exangle by Proposition 3.3 in \cite{HLN}.
\qed

\begin{proposition}{\rm \cite[Proposition 4.37]{HLN}}\label{prop2}
Let $(\C, \E, \s)$  be an $n$-exangulated category. Assume that any $\s$-inflation is
monomorphic, and any $\s$-deflation is epimorphic in $\C$. Note that this is equivalent to
assuming that any $\s$-conflation is $n$-exact sequence. We denote the class of all $\s$-conflations by $\X$. If $(\C, \E, \s)$  satisfies the following conditions {\rm (a)} and {\rm (b)} for any pair of morphisms
$A\xrightarrow{~a~}B\xrightarrow{~b~}C$ in $\C$, then $(\C, \X)$  becomes an
n-exact category in the sense of Jasso {\rm \cite[Definition 4.2]{Ja}}.

{\rm (a)} If $b\circ a$ is an $\s$-inflation, then so is $a$.

{\rm (b)} If $b\circ a$  is an $\s$-deflation, then so is $b$.

\end{proposition}

In Theorem \ref{main1}, we know that $(\A,\E_{\A},\t)$
is an $n$-exangulated category.
We assume that any $\t$-inflation is monomorphic, and any $\t$-deflation is epimorphic in $\A$. Note that this is equivalent to assuming that any $\t$-conflation is $n$-exact sequence. Thus we denote the class of all $\t$-conflations by $\mathscr{E}_{\A}$.

Our second main result is the following.

\begin{theorem}\label{main2}
Let $((\C,\E,\s)$ be a Krull-Schmidt $n$-exangulated category with enough projectives and enough injectives,
$\A$ be an $n$-extension closed subcategory of $\C$.
If $\C(\Sigma\A,\A)=\C(\A,\Omega\A)=0$, then $(\A,\mathscr{E}_{\A})$
is an $n$-exact category.
\end{theorem}
\proof By Theorem \ref{main1}, we know that $(\A,\E_{\A},\t)$
is an $n$-exangulated category.

\textbf{Step 1:} We first claim that any $\t$-inflation is monomorphic in $\A$.

Assume that $f\colon A\to B$ is a $\t$-inflation in $\A$.
By definition  there exists a distinguished $n$-exangle
$$A\xrightarrow{f}B\xrightarrow{f_1}C_2\xrightarrow{f_2}C_3\xrightarrow{f_3}\cdots
\xrightarrow{f_{n-1}}C_n\xrightarrow{f_n}C_{n+1}\overset{\delta}{\dashrightarrow}$$
with terms $C_2,C_3,\cdots,C_{n+1}\in\A$.

Now we prove that $f$ is monomorphism in $\A$.
Let $h\colon C\to A$ be a morphism in $\A$ such that $fh=0$. Since $\C$ has enough projectives, then we have the following commutative diagram of distinguished $n$-exangles by  the dual of Lemma \ref{a2}
$$\xymatrix{
\Omega C_{n+1}\ar[r]^{g_0}\ar@{}[dr] \ar@{-->}[d]^{\varphi_0} &P_1 \ar[r]^{g_1} \ar@{}[dr]\ar@{-->}[d]^{\varphi_1}&P_2\ar[r]^{g_2} \ar@{-->}[d]^{\varphi_2}\ar@{}[dr]&P_3\ar[r]^{g_3} \ar@{-->}[d]^{\varphi_3}\ar@{}[dr]&\cdot\cdot\cdot  \ar[r]^{g_{n-1}} \ar@{}[dr]& P_n\ar[r]^{g_{n}} \ar@{}[dr]\ar[d]^{\varphi_n}& C_{n+1}  \ar@{}[dr] \ar@{=}[d]^{}\ar@{-->}[r]^-{\eta} &\\
A\ar[r]^{f}&B \ar[r]^{f_1} & C_2\ar[r]^{f_2}  & C_2\ar[r]^{f_2}  & \cdot\cdot\cdot\ar[r]^{f_{n-1}}  & C_n \ar[r]^{f_{n}}&C_{n+1} \ar@{-->}[r]^-{\delta} &}$$
such that
$$M_{\bullet}:\Omega C_{n+1}\xrightarrow{\left[
                                            \begin{smallmatrix}
                                              -g_0 \\
                                              \varphi_0 \\
                                            \end{smallmatrix}
                                          \right]}
 P_1\oplus A\xrightarrow{\left[
                                            \begin{smallmatrix}
                                              -g_{1} & 0 \\
                                              \varphi_{1} & f
                                            \end{smallmatrix}
                                          \right]}
P_2\oplus B\xrightarrow{}
\cdots\xrightarrow{}
P_{n}\oplus C_{n-1}\xrightarrow{}
C_{n}\overset{}{\dashrightarrow}$$
is a distinguished $n$-exangle.

Apply the functor $\C(C,-)$ to the above distinguished $n$-exangle $M_{\bullet}$, we have the following exact sequence:
$$\C(C,\Omega C_{n+1})\xrightarrow{\C(C,\left[
                                            \begin{smallmatrix}
                                              -g_0 \\
                                              \varphi_0 \\
                                            \end{smallmatrix}
                                          \right])}\C(C,P_1\oplus A)\xrightarrow{\C(C,\hspace{0.8mm}\left[
                                            \begin{smallmatrix}
                                              -g_{1} & 0 \\
                                              \varphi_{1} & f
                                            \end{smallmatrix}
                                          \right])}\C(C,P_2\oplus B).$$
Let $\left[
                                            \begin{smallmatrix}
                                              0 \\
                                              h \\
                                            \end{smallmatrix}
                                          \right]\in\C(C,P_1\oplus A)$. Since $fh=0$, then $\left[
                                            \begin{smallmatrix}
                                              -g_{1} & 0 \\
                                              \varphi_{1} & f
                                            \end{smallmatrix}
                                          \right]\left[
                                            \begin{smallmatrix}
                                              0 \\
                                              h \\
                                            \end{smallmatrix}
                                          \right]=\left[
                                            \begin{smallmatrix}
                                              0 \\
                                              fh \\
                                            \end{smallmatrix}
                                          \right]=0$. Thus there exists a morphism $u:C \rightarrow \Omega C_{n+1}$, such that $\left[
                                            \begin{smallmatrix}
                                              -g_0 \\
                                              \varphi_0 \\
                                            \end{smallmatrix}
                                          \right]u=\left[
                                            \begin{smallmatrix}
                                             0 \\
                                              h \\
                                            \end{smallmatrix}
                                          \right]$, so $h=\varphi_0u$. That is to say, we have the following  commutative diagram

$$\xymatrix{
&  C\ar@{-->}[dl]_{u} \ar[dr]^{0}\ar[d]^{\left[
                                            \begin{smallmatrix}
                                              0 \\
                                              h \\
                                            \end{smallmatrix}
                                          \right]} & \\
\Omega C_{n+1}  \ar[r]_{\left[
                                            \begin{smallmatrix}
                                              -g_0 \\
                                              \varphi_0 \\
                                            \end{smallmatrix}
                                          \right]} &P_1\oplus A \ar[r]_{\left[
                                            \begin{smallmatrix}
                                              -g_{1} & 0 \\
                                              \varphi_{1} & f
                                            \end{smallmatrix}
                                          \right]} &P_2\oplus B
}
$$

Since $\C(\A,\Omega\A)=0$, we have $u=0$. Then $h=\varphi_0u=0$. This shows that $f$ is monomorphism in $\A$.

In a similar way, one can show that any $\t$-deflation is epimorphic in $\A$.

\textbf{Step 2:} For any pair of morphisms
$A\xrightarrow{~a~}B\xrightarrow{~b~}C$ in $\A$,
we claim that if $b\circ a$ is a $\t$-inflation, then so is $a$.
Indeed, since $h_0:=b\circ a$ is a $\t$-inflation, by definition  there exists a distinguished $n$-exangle
$$A\xrightarrow{~h_0~}C\xrightarrow{h_1}X_2\xrightarrow{h_2}\cdots\xrightarrow{h_{n-1}}X_n\xrightarrow{h_n}X_{n+1}\overset{\delta}{\dashrightarrow}$$
with terms $X_2,X_3,\cdots,X_n,X_{n+1}\in\A$.

For the $\E$-extension $a_{\ast}\del\in\E(X_{n+1},B)$ with $X_{n+1},B\in\A$,
since $\A$ is $n$-extension closed, then there exists a distinguished $n$-exangle
$$B\xrightarrow{d_0}Y_1\xrightarrow{d_1}Y_2\xrightarrow{d_2}
\cdots\xrightarrow{d_{n-1}}Y_{n}\xrightarrow{d_{n}}X_{n+1}\overset{a_{\ast}\del}{\dashrightarrow}$$
with $Y_1,Y_2,\cdots,Y_{n}\in\A$.
By (EA2$\op$) we can observe that $(a,\id_{X_{n+1}})$ has a {\it good lift} $\varphi^{\mr}=(a,\varphi_1,\ldots,\varphi_n,\id_{X_{n+1}})$, that is,
there exists the following commutative diagram of distinguished $n$-exangles
$$\xymatrix{
A \ar[r]^{h_0}\ar[d]^{a} & C \ar[r]^{h_1}\ar@{-->}[d]^{\varphi_1} & X_2 \ar[r]^{h_2} \ar@{-->}[d]^{\varphi_2}& \cdots \ar[r]^{h_{n-1}} & X_{n} \ar[r]^{h_{n}}\ar@{-->}[d]^{\varphi_n}&X_{n+1}\ar@{=}[d]\ar@{-->}[r]^-{\delta} & \\
B \ar[r]^{d_0} & Y_1 \ar[r]^{d_1} & Y_2\ar[r]^{d_2} & \cdots \ar[r]^{d_{n-1}} & Y_n \ar[r]^{d_{n}} & X_{n+1}\ar@{-->}[r]^-{a_{\ast}\del} & }$$
such that
$$A\xrightarrow{\left[
              \begin{smallmatrix}
               -h_0 \\a
              \end{smallmatrix}
            \right]}C\oplus B\xrightarrow{~}X_2\oplus Y_1\xrightarrow{~}
\cdots\xrightarrow{~}X_{n}\oplus Y_{n-1}\xrightarrow{~}Y_{n}\overset{}{\dashrightarrow}$$
is a distinguished $n$-exangle.
We observe that $X_2\oplus Y_1,\cdots,X_{n}\oplus Y_{n-1},Y_{n}\in\A$.
This shows that $\left[
              \begin{smallmatrix}
               -h_0 \\a
              \end{smallmatrix}
            \right]$ is a $\t$-inflation.
Since $\left[
              \begin{smallmatrix}
               0&1\\
                -1&0
              \end{smallmatrix}
            \right]\colon C\oplus B\to B\oplus C$ is an isomorphism, it is in particular a
$\t$-inflation. Thus
$\left[\begin{smallmatrix}
              a \\h_0
              \end{smallmatrix}
            \right]=\left[
              \begin{smallmatrix}
               0&1\\
                -1&0
              \end{smallmatrix}
            \right]\left[
              \begin{smallmatrix}
               -h_0 \\a
              \end{smallmatrix}
            \right]$  is also a $\t$-inflation.

Since $\left[
              \begin{smallmatrix}
               1&0\\
                -b&1
              \end{smallmatrix}
            \right]\colon B\oplus C\to B\oplus C$ is an isomorphism, it is in particular
a $\t$-inflation. Thus
$\left[\begin{smallmatrix}
              a \\0
              \end{smallmatrix}
            \right]=\left[
              \begin{smallmatrix}
               1&0\\
                -b&1
              \end{smallmatrix}
            \right]\left[
              \begin{smallmatrix}
               a \\h_0
              \end{smallmatrix}
            \right]$ is also a $\t$-inflation.
By definition there exists a distinguished $n$-exangle
$$A_{\bullet}:~A\xrightarrow{\left[
              \begin{smallmatrix}
               a \\0
              \end{smallmatrix}
            \right]}B\oplus C\xrightarrow{[f_1,\hspace{0.8mm} q_1]}A_2\xrightarrow{f_2}A_3\xrightarrow{f_3}\cdots\xrightarrow{f_{n-1}}A_n\xrightarrow{f_n}A_{n+1}\overset{}{\dashrightarrow}$$
with $A_2,A_3,\cdots,A_{n+1}\in\A$. By Lemma \ref{Key2}, we obtain $A_{\bullet}\simeq A'_{\bullet}\oplus C_{\bullet}$, where
$$A'_{\bullet}:~A\xrightarrow{~a~}B\xrightarrow{f_{11}}A'_2\xrightarrow{f_{21}}A_3\xrightarrow{f_3}\cdots\xrightarrow{f_{n-1}}A_n\xrightarrow{f_n}A_{n+1}\overset{}{\dashrightarrow}$$
and
$$C_{\bullet}:~0\xrightarrow{}C\xrightarrow{~1~}C\xrightarrow{}0\xrightarrow{}\cdots\xrightarrow{}0\xrightarrow{}0\xrightarrow{}0$$
is two distinguished $n$-exangles.
Since $A_2\in\A$, we also have $A'_2\in\A$.
This shows that $a$ is also a $\t$-inflation.

Similarly, we also can prove that if $b\circ a$ is a $\t$-deflation, then so is $b$.

By Proposition \ref{prop2}, we have that $(\A,\mathscr{E}_{\A})$
is an $n$-exact category.
\qed

\begin{remark}
For the proof of Theorem \ref{main2}, our proof method is to avoid proving that $\mathscr{E}_{\A}$ is closed under weak isomorphisms, which is very complicated and difficult.
\end{remark}

By applying Theorem \ref{main2} to $(n+2)$-angulated category, we have the following.

\begin{corollary}{\rm\cite[Theorem 3.2]{K}}\label{cor1}
Let $(\C,\Sigma,\Theta)$ be a Krull-Schmidt $(n+2)$-angulated category and
$\A$ be an $n$-extension closed subcategory of $\C$.
If $\C(\Sigma\A,\A)=0$, then $(\A,\mathscr{E}_{\A})$
is an $n$-exact category.
\end{corollary}

As a special case of Theorem \ref{main2}, when $n=1$, we have the following.
\begin{corollary}\label{cor}
Let $(\C,\E,\s)$ be a Krull-Schmidt extriangulated category and
$\A$ be an extension closed subcategory of $\C$.
If $\C(\Sigma\A,\A)=\C(\A,\Omega\A)=0$, then $(\A,\mathscr{E}_{\A})$
is an exact category.
\end{corollary}

\begin{remark}
In Corollary \ref{cor}, when $\C$ is a triangulated category, it is just the main theorem of \cite{D} and it was also rediscovered in \cite[Proposition 2.5]{J}.
\end{remark}

\textbf{Jian He}\\
Department of Mathematics, Nanjing University, 210093 Nanjing, Jiangsu, P. R. China\\
E-mail: \textsf{jianhe30@163.com}\\[0.3cm]
\textbf{Panyue Zhou}\\
College of Mathematics, Hunan Institute of Science and Technology, 414006 Yueyang, Hunan, P. R. China.\\
E-mail: \textsf{panyuezhou@163.com}

\end{document}